\documentclass[10pt,leqno]{article}

\usepackage{amsmath,amssymb,amsthm,mathrsfs,dsfont}

\usepackage[margin=3cm]{geometry} 

\usepackage{titlesec,hyperref}

\usepackage{color}

\usepackage{fancyhdr}
\titleformat{\subsection}{\it}{\thesubsection.\enspace}{1pt}{}

\newtheorem{theo}{Theorem}[section]
\newtheorem{lemm}[theo]{Lemma}
\newtheorem{defi}[theo]{Definition}

\newtheorem{prop}[theo]{Proposition}

\numberwithin{equation}{section}

\allowdisplaybreaks 

\begin{document}
\title{Well-posedness and global  solutions to the higher order Camassa-Holm equations with fraction-
alinertia operator in Besov space
\hspace{-4mm}
}

\author{Weikui $\mbox{Ye}^1$\footnote{E-mail:  904817751@qq.com} \quad and
\quad
 Zhaoyang $\mbox{Yin}^{2,3}$\footnote{E-mail: mcsyzy@mail.sysu.edu.cn}\\
 $^1\mbox{Institute}$ of Applied Physics and Computational Mathematics,\\
P.O. Box 8009, Beijing 100088, P. R. China\\
 $^2\mbox{Department}$ of Mathematics,
Sun Yat-sen University,\\ Guangzhou, 510275, China\\
$^3\mbox{Faculty}$ of Information Technology,\\ Macau University of Science and Technology, Macau, China}
\date{}
\maketitle
\hrule

\begin{abstract}
In this paper, we study well-posedness and  the global  solutions to the higher-order Camassa-Holm equations with fractional inertia operator in Besov space. When $a\in[\frac{1}{2},1),$ we prove the existence of the solutions in space $B^s_{p,1}(\mathbb R)$ with $s\geq 1+\frac{1}{p}$ and $p <\frac{1}{a-\frac{1}{2}}$, the existence and uniqueness of the solutions  in space $B^s_{p,1}(\mathbb R)$ with $s\geq  1+2a-\min\{\frac{1}{p},\frac{1}{p'}\},$ and the local well-posedness in space $B^s_{p,1}(\mathbb R)$ with $s>  1+2a-\min\{\frac{1}{p},\frac{1}{p'}\}$. When $a>1,$ we obtain the  existence of the solutions in space $B^s_{p,1}(\mathbb R)$ with $s\geq a+\max\{\frac{1}{p},\frac{1}{2}\}$ and the local well-posedness in space $B^s_{p,1}(\mathbb R)$ with $s\geq 1+a+\max\{\frac{1}{p},\frac{1}{2}\}$. Moreover, we obtain two results about  the global solutions. \\

\vspace*{5pt}

\vspace*{5pt}
\end{abstract}

\vspace*{10pt}

\tableofcontents
\section{Introduction and Main results}
  The two-component higher order Camassa-Holm systems with fractional inertia operator was firstly introduced
 by Escher and Lyons to describe the geometrical geodesic flow on an appropriate infinite dimensional Lie group \cite{E-Lyons}, especially for the cases that $s$ is a fractional number larger than $1$. It is the following system:
 \begin{equation}\label{system1}
\left\{\begin{array}{ll}m_{t}+2u_{x}m+um_x =\beta u_x-\kappa\rho\rho_x
 ,&t > 0,\,x\in \mathbb{R},\\
 \rho_{t}+(u\rho )_{x}=0, &t > 0,\,x\in \mathbb{R},\\
 m(t,x)=(1-\partial_x^2)^au(t,x), &t \geq 0,\,x\in \mathbb{R},\\
 u(0,x)=u_0(x), &x\in \mathbb{R},\\
 \rho(0,x)=\rho_0(x),&x\in \mathbb{R},
\end{array}\right.
\end{equation}
where $s\geq 0$ is a constant, $\beta $ is a constant which represents the vorticity of underlying flow, and $\kappa>0$ is an arbitrary real  parameter.
This model attracted a lot of  attention. In \cite{Chen-Zhou}, the authors study the local well-posedness in Besov spaces with high-regularity. In \cite{H-Y}, the authors established the local-wellposedness in Besov spaces with lower-regularity, they also give the global solutions  of \eqref{system1}  for the case $a=2,\kappa>0$. In \cite{H-Y-2}, the authors study the local and global analyticity and Gevrey regularity of \eqref{system1}.

We can take it as the generalization of the shallow water wave systems introduced in the nonlinear dynamics to describe nonlinear phenomena in  shallow water.

\par Indeed,
for $a=1$, System \eqref{system1} becomes the  two-component Camassa-Holm system:
\begin{equation}\label{2CH}
\left\{\begin{array}{ll}m_{t}+2u_{x}m+um_x =\beta  u_x-\kappa\rho\rho_x
 ,&t > 0,\,x\in \mathbb{R},\\
 \rho_{t}+(u\rho )_{x}=0, &t > 0,\,x\in \mathbb{R},\\
 m(t,x)=u(t,x)-u_{xx}(t,x), &t \geq 0,\,x\in \mathbb{R},\\
 u(0,x)=u_0(x), &x\in \mathbb{R},\\
 \rho(0,x)=\rho_0(x),&x\in \mathbb{R}.
\end{array}\right.
\end{equation}
The Cauchy problems of 2CH
 \eqref{2CH}  with
$\kappa=-1$ and with $\kappa=1$ have been discussed in \cite{E-L-Y}
and \cite{C-I}, respectively.  A new global existence
result and several new blow-up results of strong solutions for the
Cauchy problem of \eqref{2CH} with $\kappa=1$ were obtained in
\cite{G-Y}. Moreover, the authors studied the existence of global weak solutions and uniqueness of conservative weak solutions  to \eqref{2CH}  in \cite{G-Y1,L-Z}, respectively.

In this paper, we consider the case with $\rho\equiv0$ and $\beta=0$ in \eqref{system1}, then we get the higher order Camassa-Holm equation with fraction-
alinertia operator:
\begin{equation}\label{1}
\left\{\begin{array}{ll}m_{t}+2u_{x}m+um_x =0
 ,&t > 0,\,x\in \mathbb{R},\\
 m(t,x)=(1-\partial_x^2)^au(t,x), &t \geq 0,\,x\in \mathbb{R},\\
 u(0,x)=u_0(x), &x\in \mathbb{R}.
\end{array}\right.
\end{equation}

For $a=1$, Eq. \eqref{1} becomes the remarkable Camassa-Holm(CH)
equation
\begin{align}\label{CH}
u_t-u_{xxt}+2u_x(u-u_{xx})+u(u-u_{xx})_x-\alpha u_x=0,
\end{align}
modeling the unidirectional propagation of shallow water waves over
a flat bottom. Here $u(t,x)$ stands for the fluid velocity at time
$t $ in the spatial $x$ direction \cite{C-H,C-L,D-G-H,I-K, Ip, J}.
The Camassa-Holm equation is also a model for the propagation of
axially symmetric waves in hyperelastic rods \cite{C-S2, D}. It has
a bi-Hamiltonian structure \cite{C1, F-F} and is completely
integrable \cite{C-H, C3}. Also there is a geometric interpretation
of Eq.(1.1) in terms of geodesic flow on the diffeomorphism group of
the circle \cite{C-K,Ko}. One of the main features of the CH equation is that it precesses the peaked wave solution
\begin{equation*}
u(t,x)=c e^{-|x-ct|}
\end{equation*}
as the solitary weak solution
\cite{C-H-H}. They are orbitally stable and interact like solitons
\cite{B-S-S,C-S1}. The peaked traveling waves replicate a
characteristic for the waves of great height -- waves of largest
amplitude that are exact solutions of the governing equations for
water waves, cf. \cite{Ci, C-Eb, T}. Recently, it was claimed in the
paper \cite{La} that the equation might be relevant to the modeling
of tsunami, see also the discussion in \cite{C-J}.
The Cauchy problem and initial-boundary value problem for the
CH equation have been studied extensively \cite{C-Ep,
C-Em,Dan, E-Y1, E-Y2, li2,L-O, Rb, Y1}. It has been shown that this
equation is locally well-posed \cite{C-Ep, C-Em,Dan,li2,L-O, Rb} for
initial data $u_{0}\in H^{r}(\mathbb{R}),\,r>\frac{3}{2}$ and locally ill-posed \cite{G-L-M-Y} for initial data $u_{0}\in H^{\frac{3}{2}}(\mathbb{R})$. More
interestingly, it has global strong solutions \cite{Cf, C-Ep,C-Em}
and also finite time blow-up solutions \cite{Cf,C-E, C-Ep,C-Em,
C-Ez,Dan,L-O,Rb}. On the other hand, it has global weak solutions in
$ H^1(\mathbb{R}) $ \cite{B-C1,B-C2,C-Ei,C-M,X-Z}. Uniqueness of the conservative weak solutions to  Camassa-Holm equation has been proved in \cite{B-C-Z}. The advantage of
the Camassa-Holm equation in comparison with the KdV equation lies
in the fact that the Camassa-Holm equation has peaked solitons and
models wave breaking \cite{C-H-H,C-E} (by wave breaking we
understand that the wave remains bounded while its slope becomes
unbounded in finite time \cite{W}).

In this paper, we first study about the local well-posedness for Eq. \eqref{1} with initial data $u_0\in B^{s}_{p,1}(\mathbb R)$ in the Hardmard sense. If $\frac{1}{2}\leq a<1$, we get the existence of the solution  in space $B^s_{p,1}(\mathbb R)$ with $s\geq 1+\frac{1}{p}$ and $p <\frac{1}{a-\frac{1}{2}}$, the existence and uniqueness of the solutions  in space $B^s_{p,1}(\mathbb R)$ with $s\geq  1+2a-\min\{\frac{1}{p},\frac{1}{p'}\},$ and the well-posedness in space $B^s_{p,1}(\mathbb R)$ with $s>  1+2a-\min\{\frac{1}{p},\frac{1}{p'}\}$. When $a>1,$ we obtain the  existence of the solutions in space $B^s_{p,1}(\mathbb R)$ with $s\geq a+\max\{\frac{1}{p},\frac{1}{2}\}$ and the well-posedness in space $B^s_{p,1}(\mathbb R)$ with $s\geq 1+a+\max\{\frac{1}{p},\frac{1}{2}\}$. Moreover, if $m_0$ satisfies some sign condition, we can obtain two results about global solution. The proof of these local-posedness results give us a method to study the  well-posedness of system \eqref{system1}, maybe we can a  better well-posedness result than that in \cite{H-Y}. And we will study  this problem in the next paper.

The first result is about the local well-posedness with $\frac{1}{2}\leq a<1$ in critical Besov space which can be stated as follow.
 \begin{theo}\label{th1}
 Let $u_0\in B^{s}_{p,1}(\mathbb{R})$ with $\frac{1}{2}\leq a<1$. There exists a $T$ such that \\
(1) Existence: If $s\geq 1+\frac{1}{p}$ and $p <\frac{1}{a-\frac{1}{2}}$, then the system \eqref{1} have a solution $u$ belongs to
 $$E_p(T)=C([0,T);B^{s}_{p,1}(\mathbb{R}))\cap C^1([0,T);B^{s-1}_{p,1}(\mathbb{R}))).$$
(2) Uniqueness: If $s\geq  1+2a-\min\{\frac{1}{p},\frac{1}{p'}\}$, then the system \eqref{1} have a unique solution $u$ belongs to $E_p(T)$.\\
(3) Continuous dependence: If $s>  1+2a-\min\{\frac{1}{p},\frac{1}{p'}\}$, then the data-to-solution map $S_t (u_0)$ is continuous in $B^{s}_{p,1}(\mathbb{R})$.
 \end{theo}

\begin{theo}\label{th2}
 Let $u_0\in B^{s}_{p,1}(\mathbb{R})$ with $a>1$. There exists a $T$ such that \\
(1) Existence: If $s\geq a+\max\{\frac{1}{p},\frac{1}{2}\}$, then the system \eqref{1} have a solution $u$ belongs to
 $$E_p(T)=C([0,T);B^{s}_{p,1}(\mathbb{R}))\cap C^1([0,T);B^{s-1}_{p,1}(\mathbb{R}))).$$
(2) Uniqueness and continuous dependence: If $s\geq 1+a+\max\{\frac{1}{p},\frac{1}{2}\}$, then the system \eqref{1} have a unique solution $u$ belongs to $E_p(T)$.
Moreover, the data-to-solution map $S_t (u_0)$ is continuous in $B^{s}_{p,1}(\mathbb{R})$.\\
 \end{theo}

The following two results are global existence for a class of special initial data.
\begin{theo}\label{globald=3}
Let $u_0\in H^s$ with $s\geq 2a+\frac{1}{2}$. Let $u(t,x)$ be the corresponding local solution of \eqref{1} with $a>1$, if the initial data $m_0:=(1-\partial_{xx})^au_0\in L^1$ and $m_0\geq 0(or \leq0)$ , then the solution $u(t,x)$ exists globally.
\end{theo}

\begin{theo}\label{globald==3}
Let $u_0\in H^s$ with $s\geq 2a+\frac{1}{2}$. Let $u(t,x)$ be the corresponding local solution of \eqref{1} with $a>1$, if the initial data $m_0:=(1-\partial_{xx})^au_0\in L^1$ and $m_0(x)$ is an odd function such that $m_0\leq 0$ when $x\leq 0$, and $m_0\geq 0$ when $x\geq 0$, then the solution $u(t,x)$ exists globally.
\end{theo}

The remainder of the paper is organized as follows. In Section 2 we introduce some preliminaries which will be used in sequel. In section 3, we prove the local-posedness of Eq. eqref{1}, i.e. Theorems 1.1-1.2. Finally, we give the proof of the existence of global solutions, i.e. Theorems 1.3-1.4.

\section{Preliminaries}
In this section, we will recall some propositions and lemmas on the Littlewood-Paley decomposition and Besov spaces.

\begin{prop}\cite{book}
Let $\mathcal{C}$ be the annulus $\{\xi\in\mathbb{R}^d:\frac 3 4\leq|\xi|\leq\frac 8 3\}$. There exist radial functions $\chi$ and $\varphi$, valued in the interval $[0,1]$, belonging respectively to $\mathcal{D}(B(0,\frac 4 3))$ and $\mathcal{D}(\mathcal{C})$, and such that
$$ \forall\xi\in\mathbb{R}^d,\ \chi(\xi)+\sum_{j\geq 0}\varphi(2^{-j}\xi)=1, $$
$$ \forall\xi\in\mathbb{R}^d\backslash\{0\},\ \sum_{j\in\mathbb{Z}}\varphi(2^{-j}\xi)=1, $$
$$ |j-j'|\geq 2\Rightarrow\mathrm{Supp}\ \varphi(2^{-j}\cdot)\cap \mathrm{Supp}\ \varphi(2^{-j'}\cdot)=\emptyset, $$
$$ j\geq 1\Rightarrow\mathrm{Supp}\ \chi(\cdot)\cap \mathrm{Supp}\ \varphi(2^{-j}\cdot)=\emptyset. $$
The set $\widetilde{\mathcal{C}}=B(0,\frac 2 3)+\widetilde{\mathcal{C}}$ is an annulus, and we have
$$ |j-j'|\geq 5\Rightarrow 2^{j}\mathcal{C}\cap 2^{j'}\widetilde{\mathcal{C}}=\emptyset. $$
Further, we have
$$ \forall\xi\in\mathbb{R}^d,\ \frac 1 2\leq\chi^2(\xi)+\sum_{j\geq 0}\varphi^2(2^{-j}\xi)\leq 1, $$
$$ \forall\xi\in\mathbb{R}^d\backslash\{0\},\ \frac 1 2\leq\sum_{j\in\mathbb{Z}}\varphi^2(2^{-j}\xi)\leq 1. $$
\end{prop}

\begin{defi}\cite{book}
Denote $\mathcal{F}$ by the Fourier transform and $\mathcal{F}^{-1}$ by its inverse.
Let $u$ be a tempered distribution in $\mathcal{S}'(\mathbb{R}^d)$. For all $j\in\mathbb{Z}$, define
$$
\Delta_j u=0\,\ \text{if}\,\ j\leq -2,\quad
\Delta_{-1} u=\mathcal{F}^{-1}(\chi\mathcal{F}u),\quad
\Delta_j u=\mathcal{F}^{-1}(\varphi(2^{-j}\cdot)\mathcal{F}u)\,\ \text{if}\,\ j\geq 0,\quad
S_j u=\sum_{j'<j}\Delta_{j'}u.
$$
Then the Littlewood-Paley decomposition is given as follows:
$$ u=\sum_{j\in\mathbb{Z}}\Delta_j u \quad \text{in}\ \mathcal{S}'(\mathbb{R}^d). $$

Let $s\in\mathbb{R},\ 1\leq p,r\leq\infty.$ The nonhomogeneous Besov space $B^s_{p,r}(\mathbb{R}^d)$ is defined by
$$ B^s_{p,r}=B^s_{p,r}(\mathbb{R}^d)=\{u\in S'(\mathbb{R}^d):\|u\|_{B^s_{p,r}(\mathbb{R}^d)}=\Big\|(2^{js}\|\Delta_j u\|_{L^p})_j \Big\|_{l^r(\mathbb{Z})}<\infty\}. $$
\end{defi}

\begin{prop}\cite{book}
(1) For any $p\in[1,\infty]$, the Besov space $B^{\frac{d}{p}}_{p,1}$ is a Banach algebra. If $f,g\in B^{\frac{d}{p}}_{p,1}$, then
$$\|fg\|_{B^{\frac{d}{p}}_{p,1}}\leq C\|f\|_{B^{\frac{d}{p}}_{p,1}}\|g\|_{B^{\frac{d}{p}}_{p,1}}.$$
(2) For any $s\in\mathbb{R}$, $p,r\in[1,\infty]$, if $u\in B^s_{p,r}$, then $\nabla u\in B^{s-1}_{p,r}$ and $(1-\Delta )^{-1} u\in B^{s+2}_{p,r}$. Moreover,
$$\|\nabla u\|_{B^{s-1}_{p,r}}\leq C\|u\|_{B^s_{p,r}}, \quad \|(1-\Delta )^{-1} u\|_{B^{s+2}_{p,r}}\leq C\|u\|_{B^s_{p,r}} $$
\end{prop}

In order to prove our main theorem, we have to use the following result about the transport equation
\begin{equation}\label{s1}
\left\{\begin{array}{l}
    f_t+v\cdot\nabla f=g, \\
    f(0,x)=f_0(x),
\end{array}\right.
\end{equation}
\begin{lemm}\label{priori estimate}\cite{book}
Let $s\in [\max\{-\frac{d}{p},-\frac{d}{p'}\},\frac{d}{p}+1]$ ($s=1+\frac{1}{p},r=1$; $s=\max\{-\frac{d}{p},-\frac{d}{p'}\},r=\infty$).
There exists a constant $C$ such that for all solutions $f\in L^{\infty}([0,T];B^s_{p,r})$ of \eqref{s1} with initial data $f_0$ in $\dot{B}^s_{p,r}$, and $g$ in $L^1([0,T];B^s_{p,r})$, we have, for a.e. $t\in[0,T]$,
$$ \|f(t)\|_{B^s_{p,r}}\leq e^{C_2 V(t)}\Big(\|f_0\|_{B^s_{p,r}}+\int_0^t e^{-C_2 V(t')}\|g(t')\|_{B^s_{p,r}}dt'\Big), $$
where $V'(t)=\|\nabla v\|_{B^{\frac{d}{p}}_{p,r}\cap L^{\infty}}$(if $s=1+\frac{1}{p},r=1$, $V'(t)=\|\nabla v\|_{B^{\frac{d}{p}}_{p,1}}$).
\end{lemm}

\begin{lemm}\label{Convergence}\cite{book}
Let $1\leq p<\infty$. Define $\overline{\mathbb{N}}=\mathbb{N}\cup \{\infty\}$. Suppose $f\in L^1([0,T];B^{\frac{d}{p}}_{p,1})$ and $a_0\in B^{\frac{d}{p}}_{p,1}$. For $n\in \overline{\mathbb{N}}$, denote by $a^n\in C([0,T];B^{\frac{d}{p}}_{p,1})$ the solution of
\begin{align}\label{an}
\left\{
\begin{array}{ll}
\partial_t a^n+A^n\cdot\nabla a^n=f,\\[1ex]
a^n|_{t=0}(x)=a_{0}(x).\\[1ex]
\end{array}
\right.
\end{align}
Assume that $\sup_{n\in\overline{\mathbb{N}} }\|A^n\|_{B^{1+\frac{d}{p}}_{p,1}}\leq \alpha(t)$ with $\alpha(t)\in L^1(0,T)$. If $A^n\to A^\infty$ in $L^1(B^{\frac{d}{p}}_{p,1})$, then
the sequence $a^n\to a^\infty$ in $C([0,T];B^{\frac{d}{p}}_{p,1})$.
\end{lemm}

{\bf Notation.} For simplicity, we drop $\mathbb{R}$ or $\mathbb{T}$ in the notation of function spaces if there is no ambiguity.
\section{Local well-posedness.}
In this section, we will prove the local well-posedness for the system \eqref{1}.

\subsection{Uniformly bound and existence.}
{\bf The proof of Theorem \ref{th1}:} Starting with $u^0\doteq 0$, we construct a sequence of the approximate solutions $u^n$ by solving the following linear transport equations:
\begin{align}\label{2n}
\left\{
\begin{array}{ll}
u^{n+1}_{t}+u^n u^{n+1}_x=[(1-\partial_{xx} )^{-a},u^n\partial_x]m^n-2(1-\partial_{xx} )^{-a}(u^n_xm^n),\\[1ex]
m^n=(1-\partial_{xx} )^{-a}u^n,\\[1ex]
u^n|_{t=0}(x)=S_nu_{0}(x).\\[1ex]
\end{array}
\right.
\end{align}

Without loss of generality, we consider the low regularity case $s=1+\frac{d}{p}$.
Assume that $u^n\in E_p(T)=C([0,T];B^{1+\frac{1}{p}}_{p,1})\cap C^1([0,T];B^{\frac{1}{p}}_{p,1})$. Taking advantage of the standard theory of the transport equations, we can deduce that there exist a unique solutions $u^{n+1} \in E_p(T)$. Using the lemma \ref{priori estimate}, we have
\begin{align}\label{3.7}
\|u^{n+1}(t)\|_{B^{1+\frac{1}{p}}_{p,1}}&\leq e^{C\int^t_0\|u^n\|_{B^{1+\frac{1}{p}}_{p,1}}ds}\|u_0\|_{B^{s}_{p,1}}\notag\\
&\quad+e^{C\int^t_0\|u^n\|_{B^{1+\frac{1}{p}}_{p,1}}ds}\int^t_0e^{-C\int^{t'}_0\|u^n\|_{B^{1+\frac{1}{p}}_{p,1}}ds}\|[(1-\partial_{xx} )^{-a},u^n\partial_x]m^n-2(1-\partial_{xx} )^{-a}(u^n_xm^n)\|_{B^{1+\frac{1}{p}}_{p,1}}dt' ,\notag\\
&\leq e^{C\int^t_0\|u^n\|_{B^{1+\frac{1}{p}}_{p,1}}ds}\bigg(\|u_0\|_{B^{s}_{p,1}}+\int^t_0e^{-C\int^{t'}_0\|u^n\|_{B^{1+\frac{1}{p}}_{p,1}}ds}\|u^n\|^2_{B^{1+\frac{1}{p}}_{p,1}}dt'\bigg),
\end{align}
where we use the fact that
$$[(1-\partial_{xx} )^{-a},u^n\partial_x]m^n=[(1-\partial_{xx} )^{-a},T_{u^n}\partial_x]m^n+(1-\partial_{xx} )^{-a}T_{m^{n}_x}u^n+T_{u^{n}_x}u^n+(1-\partial_{xx} )^{-a}R(u^n,m^{n}_x)+R(u^n,u^{n}_x),$$
by Bony decomposition and Lemma 10.25 in \cite{book} we get the above estimation.

This implies that there exist a $T$ satisfies that $0<T<\frac{1}{C\|u_0\|^2_{B^{1+\frac{1}{p}}_{p,1}}}$ such that
\begin{align}\label{3.8}
\sup_{t\in[0,T]}\|u^{n+1}\|_{B^{1+\frac{1}{p}}_{p,1}}\leq \frac{C\|u_0\|_{B^{1+\frac{1}{p}}_{p,1}}}{1-CT\|u_0\|^2_{B^{1+\frac{1}{p}}_{p,1}}}.
\end{align}
Therefore, $\{u^n\}_{n\in\mathbb{N}}$ is uniformly bounded in $L^\infty_T(B^{1+\frac{d}{p}}_{p,1})$. From the system \eqref{2n}, we can deduce that $\partial_t u^n$ is uniformly bounded in $L^\infty_T(B^{\frac{1}{p}}_{p,1})$. Using an interpolation argument, we obtain that $u^n$ is uniformly bounded in $C([0,T);B^{1+\frac{1}{p}}_{p,1})\cap C^{\frac{1}{2}}([0,T);B^{\frac{1}{p}}_{p,1})$. Taking advantage of Cantor's diagonal process and Ascoli's theorem, we can obtain a function $u_j$ such that $\phi_ju^n$ converges to $u_j$ with $\phi_j$ is a smooth function with support in the ball $B(0,j+1)$. Moreover, we can verify that there exists a function $u$ such that for all $\phi\in C^\infty_0$, $\phi u^n\rightarrow \phi u$.(For more details, see the Chapter 10 in \cite{book}). By the Fatou property and the interpolation argument, one can check that $u\in E_p(T)$ is the solution of the system \eqref{1}.

\subsection{Uniqueness.}
We will prove the uniqueness of solutions to \eqref{1} next. Actually, to prove the uniqueness we need more regularity such as $u\in C_T(B^{2+\frac{1}{p}}_{p,1})$.
let's recall \eqref{1}:
 \begin{align}\label{2}
\left\{
\begin{array}{ll}
m_{t}+um_x+2u_x m=0,\\[1ex]
u=(1-\partial_{xx})^{-a}m,\quad a\in\mathbb{R}^+/\mathbb{Z}^+\\[1ex]
m|_{t=0}(x)=m_{0}(x),\\[1ex]
\end{array}
\right.
\end{align}
Since the form of \eqref{2} is more simple than \eqref{1}, one can estimate on \eqref{2}. This need less regularity since $2+\frac{1}{p}\leq 1+2a-\min\{\frac{1}{p},\frac{1}{p'}\},~\frac{1}{2}\leq a<1$. Suppose that $m_i=(1-\partial_{xx})^{-a}u_i,i=1,2$ are two solutions of \eqref{2}, then $u_i\in C_T(B^{1+2a-\min\{\frac{1}{p},\frac{1}{p'}\}}_{p,1}))~(i=1,2)$ are two solutions of \eqref{2}. Setting $w=m_1-m_2$, we obtain
\begin{equation*}
    \partial_tw+m_1\partial_xw=-wm_{2x}+2u_{1x}w+2(1-\partial_{xx})^{-a}w_xm_2,
\end{equation*}

By virtue to Lemma \ref{priori estimate} and the Bony decomposition, we have
\begin{align}\label{ineq11}
    \|w(t)\|_{B^{-\min\{\frac{1}{p},\frac{1}{p'}\}}_{p,\infty}}&\leq C_{u_0}(\|w(0)\|_{B^{-\min\{\frac{1}{p},\frac{1}{p'}\}}_{p,\infty}}
    +\int_{0}^{t}\|-wm_{2x}+2u_{1x}w+2(1-\partial_{xx})^{-a}w_xm_2\|_{B^{-\min\{\frac{1}{p},\frac{1}{p'}\}}_{p,\infty}}ds)\notag\\
    &\leq C_{u_0}(\|w(0)\|_{B^{-\min\{\frac{1}{p},\frac{1}{p'}\}}_{p,\infty}}
    +\int_{0}^{t}\|u_1,u_2\|_{B^{1+2a-\min\{\frac{1}{p},\frac{1}{p'}\}}_{p,1}}\|w\|_{B^{-\min\{\frac{1}{p},\frac{1}{p'}\}}_{p,\infty}}ds)\notag\\
    &\leq C_{u_0}(\|w(0)\|_{B^{-\min\{\frac{1}{p},\frac{1}{p'}\}}_{p,\infty}})
\end{align}

Therefore, the uniqueness is obvious in view of \eqref{2}.

\subsection{Continuous dependent.}
This subsection devote to study about the continuous dependent on initial.
Assume that $u^n_0\to u^\infty_0$ in $B^{s}_{p,1},s>1+2a-\min\{\frac{1}{p},\frac{1}{p'}\}$ and $u^n,u^\infty$ are the solutions of \eqref{1} with the initial data $u^n_0,u^\infty_0$ respectively. Notice that their corresponding solutions $u^n, u^\infty$ are uniformly bounded in $L^\infty_T(B^{1+\frac{d}{p}}_{p,1})$. An interpolation
argument and \eqref{ineq11} yield that $u^n\to u^\infty$ in $C([0,T);B^{s-\varepsilon}_{p,1})$ for any $\varepsilon>0$.
In order to prove that $u^n \to u^\infty$ in $C([0,T);B^{s}_{p,1})$, it sufficient to show that $m^n_x \to   m^\infty_x$ in $C([0,T);B^{s-1-2a}_{p,1})$.
Let $V^n=m^n_x$ for all $n\in \overline{\mathbb{N}}$. Split $V^n$ into $W^n+Z^n$ with $(W^n,Z^n)$ satisfying that
  \begin{align*}
\left\{
\begin{array}{ll}
W^n_{t}+u^nW^{n}_x=-u^{\infty}_xm^{\infty}_x+\partial_x{u^{\infty}_xm^{\infty}}:=F^{\infty},\\[1ex]
W^n|_{t=0}(x)=V^n_0=m^{\infty}_{x0}(x).\\[1ex]
\end{array}
\right.
\end{align*}
and
  \begin{align*}
\left\{
\begin{array}{ll}
Z^n_{t}+u^nZ^{n}_x=F^n-F^{\infty} ,\\[1ex]
Z^n|_{t=0}(x)=m^n_{x0}- m^{\infty}_{x0}.\\[1ex]
\end{array}
\right.
\end{align*}
By virtue of the lemma \eqref{Convergence} and $(1-\partial_{xx})^{-a}u^n_x=m^n_x=V^n=W^n+Z^n$, we verify that
\begin{align}
W^n \to W^\infty   \quad \quad \text{in} \quad \quad  C([0,T];B^{s-1-2a}_{p,1}).
\end{align}
Using Bony decomposition, we deduce that
\begin{align}
\|F^n-F^{\infty}\|_{B^{s-1-2a}_{p,1}}&\leq C(\|W^n-W^{\infty}\|_{B^{s-1-2a}_{p,1}}+\|z^n-0\|_{B^{s-1-2a}_{p,1}})\|m^n,m^{\infty}\|_{B^{s-2a}_{p,1}}\notag\\
&\leq C_{u_0}(\|W^n-W^{\infty}\|_{B^{s-1-2a}_{p,1}}+\|z^n-0\|_{B^{s-1-2a}_{p,1}}).
\end{align}
Indeed, the uniqueness of transport equation ensures that $z^\infty \equiv 0$. Taking advantage of lemma \eqref{priori estimate}, we get
\begin{align}
\begin{split}
\|z^n\|_{B^{s-1-2a}_{p,1}}&\leq C(\|u^n_0-u^\infty_0\|_{B^{1+\frac{d}{p}}_{p,1}}+\int^t_0 \|F^n-F^{\infty}\|_{B^{s-1-2a}_{p,1}}  d\tau)\\
&\leq C_{u_0}(\|u^n_0-u^\infty_0\|_{B^{1+\frac{d}{p}}_{p,1}}+\int^t_0\|W^n-W^\infty\|_{B^{s-1-2a}_{p,1}}+\|z^n-0\|_{B^{s-1-2a}_{p,1}} d\tau) .
\end{split}
\end{align}
Using the facts that $\lim_{n\to \infty}\|u^n_0-u^\infty_0\|_{B^{s}_{p,1}}=0 $, $\lim_{n\to \infty} \|u^n-u^\infty\|_{B^{s-1}_{p,1}}=0$, $\lim_{n\to \infty }\|W^n-W^\infty\|_{B^{s-1-2a}_{p,1}}=0$, and the Gronwall inequality yields that $z_n$ tends to $0$ in $C([0,T];B^{s-1-2a}_{p,1})$. Since $(1-\partial_{xx})^{-a}u^n_x=m^n_x=V^n=W^n+Z^n$, it follows that
\begin{align}
\begin{split}
\|V^n-V^\infty \|_{B^{\frac{d}{p}}_{p,1}}&\leq C(\|W^n-W^\infty\|_{B^{s-1-2a}_{p,1}}+\|Z^n-Z^\infty\|_{B^{s-1-2a}_{p,1}})\\
&\leq  C (\|W^n-W^\infty\|_{B^{s-1-2a}_{p,1}}+\|Z^n\|_{B^{s-1-2a}_{p,1}}) \rightarrow 0  \quad \text{as} \quad n\to \infty,
\end{split}
\end{align}
that is $m^n_x \to m^\infty_x$ in $C([0,T);B^{s-1-2a}_{p,1})$. Therefore we complete the proof of Theorem \ref{th1}. Similarly, we can obtain Theorem \ref{th2}.

\section{Global existence.}
In this section we construct a class of special data such that the corresponding solution is global in time when $a>1$.

{\bf Proof of theorem \ref{globald=3}:}
\begin{proof}
Set $T$ be the maximal time of $m(t,x)$.
Firstly,
one can easily deduce that $m(t,x)\geq 0$ if $m_0(x)\geq 0$ by the characteristic method. Since the functional calculus tells us that
$$u(t,x)=(1-\partial_{xx})^{-a}m(t,x)=C\int_{0}^{+\infty}s^{a-1}e^{-s}e^{s\Delta}m(t,x)ds:=(G*m)(t,x),$$
we deduce that $u(t,x)\geq 0$ if $m_0(x)\geq 0$. Moreover, we have $G(z)=G(-z)$.

Then, by $m_0\in L^1(\mathbb{R})$, one have
\begin{align}
\frac{d}{dt}\|m(t)\|_{L^1}=\frac{d}{dt}\int_{\mathbb{R}}m(t,x)dx&=-\int_{\mathbb{R}}\partial_x (um)(t,x)dx-\int_{\mathbb{R}}(u_xm)(t,x)dx\notag\\
&=0-<G*m_x,m>\notag\\
&=0+<G*m,m_x>\notag\\
&=0+<m,G*m_x>,
\end{align}
where the last inequality holds by $G(x)=G(-x)$. So we have $-<G*m_x,m>=+<m,G*m_x>=0$ and
$$\frac{d}{dt}\|m(t)\|_{L^1}=0.$$

Finally, since
$$\|u_x\|_{L^{\infty}}\leq C \|m\|_{B^{1-2a}_{\infty,1}}\leq C \|m\|_{B^{2-2a}_{1,1}}\leq C \|m\|_{L^1}\leq C \|m_0\|_{L^1},~~if~~a>1, $$
we obtain the global existence by the blow-up criteria and the bootstrap argument.

\end{proof}
Next we prove another global solution of \eqref{1} with different form of the initial data.\\
{\bf Proof of theorem \ref{globald==3}:}
\begin{proof}
Set $T$ be the maximal time of $m(t,x)$.
Firstly, since $u(t,x)=(G*m)(t,x)$ and $G(z)=G(-z)$, one can easily deduce that $m(t,x)$ is an odd function if $m_0(x)$ is an odd function.
Similar to the proof of Theorem \ref{globald=3}, by the characteristic method we deduce that if $m_0\leq 0$ when $x\leq 0$, $m_0\geq 0$ when $x\geq 0$, then we obtain
\begin{align}\label{g1}
m_(t,x)\leq 0~~when~~ x\leq q(t,0);~~ m_0\geq 0~~ when~~ x\geq q(t,0),~~\forall t\in [0,T),
\end{align}
where $q(t,\xi)=\xi+\int_{0}^{t}m(s,q(s,\xi))ds$ is the characteristic curves.

Next we want to prove that
\begin{align}\label{g2}
m(t,x)\leq 0~~when~~ x\leq 0,~~ m_0\geq 0~~ when~~ x\geq 0,~~\forall t\in [0,T).
\end{align}
If $q(t,0)=0$ for any $t\in[0.T)$, then by \eqref{g1} we immediately get \eqref{g2}. Other if $q(t,0)>0$ (or $q(t,0)<0$) for some $t\in[0.T)$, by \eqref{g1} we deduce that
$$m(t,x)\leq 0~~when~~ x\in [-q(t,0),0];~~ m_0\leq 0~~ when~~ x\in [0,q(t,0)].$$
Since $m(t,x)$ is an odd function such that $m(t,x)=-m(t,-x),~x\in [-q(t,0),q(t,0)]$, we obtain
$$m(t,x)= 0~~when~~ x\in [-q(t,0),q(t,0)].$$
By \eqref{g1} again we still get \eqref{g2}.

Finally, by \eqref{g2} we can easily deduce that
$\frac{d}{dt}\|m(t)\|_{L^1}=0$ and $\|u_x\|_{L^{\infty}}\leq C \|m_0\|_{L^1}$,~~for~~$a>1$. By the blow-up criteria and the bootstrap argument we obtain the global existence.

\end{proof}

 \noindent\textbf{Acknowledgements.} This work was partially supported by National Natural Science Foundation
of China [grant number 11671407 and 11701586], the Macao Science and Technology Development Fund (grant number 0091/2018/A3), Guangdong Special Support Program (grant number 8-2015), and the key project of NSF of  Guangdong province (grant number 2016A030311004).

\bibliographystyle{abbrv} 

\end{document}